\documentclass[12pt]{article}

\newtheorem{Theorem}{Theorem}            
\newtheorem{Proposition}{Proposition}[section]  
\newtheorem{Lemma}[Proposition]{Lemma}
\newtheorem{Corollary}[Proposition]{Corollary}

\newtheorem{Statement}[Proposition]{Statement}

\def\C{{\mathcal C}}

\newcommand{\be}{\begin{eqnarray*}}
\newcommand{\ee}{\end{eqnarray*}}
\newcommand{\beeq}{\begin{equation}}
\newcommand{\eeeq}{\end{equation}}
\newcommand{\bee}{\begin{eqnarray}}
\newcommand{\eee}{\end{eqnarray}}

\newcommand{\ve}{\varepsilon}

\newcommand{\xa}{\xi_{\alpha}}

\newcommand{\prob} {(\Omega, {\mathcal F}, P)}

\def	\O	{{\Omega}}
\def	\o      {{\omega}}

\def\qed{\hbox{\hskip 6pt\vrule width6pt height7pt depth1pt
    \hskip1pt}\bigskip}

\def\bbbr{{\rm I\!R}}
\newcommand{\R}{\bbbr}

\begin{document}

\begin{center}{{\bf\Large Stieltjes integrals of H\"older continuous functions 
with applications to fractional Brownian motion.} }
\vspace{5mm}

{\large A.\ A.\ Ruzmaikina}\footnote{ Present address: Department of Mathematics, University of Virginia, Kerchof Hall, Charlottesville VA 22903. Email: anastasia@virginia.edu.}
\vspace{5mm}

\begin{center}
{\small\it
Department of Physics, Princeton University, Princeton NJ 08544  \\
}
\end{center}

\end{center}






\begin{abstract}
We give a new estimate on Stieltjes integrals of H\"older continuous functions and use it to prove an existence-uniqueness theorem for solutions of ordinary differential equations with H\"older continuous forcing.
 We construct stochastic integrals with respect to fractional 
Brownian motion, and establish sufficient conditions for its
existence. 
We prove that stochastic differential equations with fractional
Brownian motion have a unique solution with probability $1$ in certain
classes of H\"older-continuous functions. We give tail estimates
of the maximum of stochastic integrals from tail estimates of
the H\"older coefficient of fractional Brownian motion. In addition we
apply the techniques used for ordinary Brownian motion to construct
stochastic integrals of deterministic functions with respect to
fractional Brownian motion and give tail estimates of its maximum.
\end{abstract}


\section{Introduction.}

	{\it Fractional Brownian motion\/} (fBm) was first introduced by
Kolmogorov in 1940 \cite{Kl} and later studied by L\'evy and
Mandelbrot \cite{Lv, Mn}.	

Let $\prob$ be a probability space, and $\alpha \in \R$, $|\alpha| < 1$ be a parameter.  FBm with exponent $\alpha$ is a self-similar, centered Gaussian random process 
$\xa(t, \o)$, $(t, \o)\in [0, \infty) \times\O$ (often abbreviated as 
$\xa(t)$) with stationary increments and the correlation function
$$
E(\xa(s),\xa(t))
=C(s^{1+\alpha}+t^{1+\alpha}-|s-t|^{1+\alpha}), \quad 
C = - \frac{\Gamma(1-\alpha)}{\alpha}\
     \frac{\cos\frac{1+\alpha}{2}\pi}{\frac{1+\alpha}{2}\pi}.
$$
For $\alpha=0$ we recover the ordinary Brownian motion (oBm).  
From the form of the correlation function of the increments
$$
E(d\xa(s),d\xa(t)) = C \alpha(\alpha + 1) \frac{ds dt}{|s - t|^{1-\alpha}},
$$
one can see that fBm is not Markovian for all $\alpha \ne 0$.  The trajectories of fBm are almost surely H\"older-continuous with exponent 
less than $\frac{1+\alpha}{2}$, and not H\"older-continuous with exponent greater or equal to
$\frac{1+\alpha}{2}$.  The finite-dimensional distributions of fBm 
are scale-invariant:
$$
\forall r, t > 0, \quad
{\textstyle  r^{- \frac{1+\alpha}{2}}\xa(rt) \stackrel{\mathrm{dist}}{=}  \xa(t)}.
$$
The scale-invariance and long-range correlations make fBm important 
in applications.

In recent years there has been much interest in fBm, see for example 
\cite{G, K1, K2, Lb, Ml, Sn, T1, T2, W, X, Yr}.
The problem of the construction of stochastic calculus with respect to fBm has been considered in \cite{Ln, DH, DU}.  
The main difficulty is
that fBm fails to be a semi-martingale for all $0 < \alpha < 1$ 
\cite{Ln}. Our goal is to construct the stochastic integrals with respect to fBm and 
to prove the general existence-uniqueness theorem for solutions of stochastic 
differential equations (SDE's) with fBm. Due to the strongly non-Markovian nature 
of fBm, the most natural way to define the stochastic integral is to do it pathwise, 
for a.e. $\omega$.  
This brings us to the question of existence of Stieltjes integrals with respect 
to H\"older continuous functions and of the existence and uniqueness of solutions 
of ODE's with H\"older continuous forcing (Section ~\ref{maintheorems}).  
The existence of Stieltjes integrals with respect to H\"older continuous functions 
was established in \cite{Yg, Kn}. Below we use the methods of Renormalization group 
to prove a formula (Theorem~\ref{sum}) which allows us to estimate the $L^{\infty}$-norm of the Stieltjes 
integral. 
We then use this estimate to prove the general existence-uniqueness theorem for 
solutions of ODE's with H\"older continuous forcing. 
In section ~\ref{applications} we apply these results to construct the  stochastic 
integrals with respect to fBm and show the
existence and uniqueness of solutions of stochastic differential 
equations with fBm. 
Previous results in this direction were obtained by \cite{Ln, DH}. 
In section~\ref{additional}  we use Theorem~\ref{sum} to get estimates on the tails of the stochastic integral with respect to fBm. 
We also consider the question of existence of fBm stochastic integral of deterministic functions and derive the probability distribution of its maximum.                                             
Throughout the paper, $\C^{\beta}(I)$ denotes the space of 
H\"older-continuous functions on the interval $I$ with exponent 
$\beta$.  We will have many occasions to partition an interval into 
$2^n$ sub-intervals of equal size; the $i$-th partition point of the 
interval under discussion is denoted by $s_i^n$, and we put 
$\Delta f(s+s_i^n) = f(s+s_{i+1}^n) - f(s+s_i^n)$ if $f$ is a function 
defined on the interval.

\section{Main Theorems.}\label{maintheorems}
In this section we consider the question of existence of a Stieltjes integral for functions of unbounded variation, give a formula which allows to estimate its upper bound and prove the existence-uniqueness theorem for ordinary differential equations with H\"older continuous forcing.

In applications (such as the construction of the stochastic calculus for fBm), it is often interesting to consider Stieltjes integrals $\int f dg$ for functions of unbounded variation. The difficulty in constructing the integral is that the upper bounds on Riemann sums $\sum |f||\Delta g|$ diverge. However this problem can be solved on certain classes of $f,g,$ since, if 
$g$ oscillates fast enough, the nearby terms in the Riemann sum $\sum f \Delta g$ may cancel. It is shown in \cite{Yg, Kn} that the Stieltjes integral exists on certain classes of H\"older continuous functions.

\begin{Theorem}\label{Stieltjes} (Young-Kondurar)
Let $f \in \C^{\beta}(\R)$, $g \in \C^{\gamma}(\R)$. If $\beta + \gamma > 1,$
then $\int_{0}^{t}fdg$ exists as a Stieltjes integral for all $t>0.$

\end{Theorem}

Generalizations of Theorem ~\ref{Stieltjes} can be found in \cite{Dy}.

The next formula is useful in estimating the $L^{\infty}$-norm
of the Stieltjes integral and is the main tool used in this paper.
\begin{Theorem}\label{sum}
Let $f$, $g$,
$\beta$, and  $\gamma$  be as in Statement~\ref{Stieltjes}.
Then 
\begin{equation}\label{theorem-sum}
\int_s^t f(\tau)dg(\tau) 
	= f(s)(g(t)-g(s)) 		
	+
	\sum_{k=1}^{\infty}
	\sum_{i=0}^{2^{k-1}-1}
	\Delta f(s+s_{2i}^k) \Delta g(s+s_{2i+1}^k).
\end{equation}
\end{Theorem}

\noindent
As far as we know this form of the Stieltjes integral has not appeared before.
The idea behind the proof of Theorem ~\ref{sum} is to write a recursion between the Riemann sums on finer partitions of the interval and the Riemann sum on coarser partitions of the interval, very much like in Renormalization group. The same idea is used in \cite{Ru} to give a new proof of Theorem ~\ref{Stieltjes}, which amounts to showing that the right hand side of the equation ~(\ref{theorem-sum}) does not depend on the sequence of partitions we choose. Below we will need the change of variables formula which is a deterministic analog of the It\^o's formula for Brownian motion.

\begin{Lemma}\label{ito}
Let $u:[0,\infty)\times \R \rightarrow \R$, $e, f:[0,\infty) \rightarrow \R$,  $1 > \gamma > \frac{1}{2}$, $\beta> 1-\gamma$, $T > 0$.
Suppose $f \in C^{\beta}([0,T])$, $g \in C^{\gamma}([0,T])$, $e$ is continuous, $u$ is differentiable in $t$ with continuous $\partial u/\partial t$ and twice differentiable in $x$.
  For $0 \leq t \leq T$, consider the Stieltjes integral
\be
	\eta(t)=\eta(0)+\int_{0}^{t}e(s)ds+\int_{0}^{t}f(s)dg(s).
\ee
Then for all $0 \leq t \leq T$, $v(t)=u(t,\eta(t))$ is also a Stieltjes integral whose differential is
\begin{equation}\label{ito1}
	dv = \frac{\partial u}{\partial t} dt + 
			\frac{\partial u}{\partial \eta} d\eta.
\end{equation}
\end{Lemma}
Note that for $\gamma > \frac{1}{2}$,  the change of variable formula is the same as in the case of ordinary calculus. This follows from the fact that the quadratic variation of $g$ is zero and therefore the terms of order $d \eta^2$ are negligible.

Theorems ~\ref{Stieltjes} and ~\ref{sum} can be used to prove the following general theorem on the existence and uniqueness of solutions of ordinary differential equations with H\"older continuous forcing.
\begin{Theorem}\label{main}
Let $b, \sigma : [0, \infty) \times \R \rightarrow \R$, $g \in \C^{\gamma}(\R)$ and $1/2 < \gamma \leq 1$. 
Suppose $b$ is globally Lipschitz in $t$ and $x$, and $\sigma \in \C^{1}(\R)$ with $\sigma$, $\sigma'_t$, $\sigma'_x$ globally Lipschitz in $t$ and $x$.
Then for every $T > 0$ and $\gamma > \beta > 1-\gamma$, the ordinary differential equation
\begin{equation}\label{sde}
	dx(t) = b(t, x(t))dt
	+ \sigma(t, x(t))d g(t), \quad x(0) = x_{0}  
\end{equation}
has a unique solution in $C^{\beta}([0,T])$.
\end{Theorem}

\section{Applications to fractional Brownian motion.}\label{applications}
The natural way of constructing the fBm stochastic integral $\int_{0}^{t}f(\tau, \o)d\xa(\tau, \o)$ is to define it as a Stieltjes integral for a.e.\ $\o$.  Since $\xa(\tau, \o) \in \C^{\gamma}(\R)$ for $\gamma \geq \frac{1+\alpha}{2}$ with probability one, Theorem ~\ref{Stieltjes} implies that the fBm stochastic integral $\int_{0}^{t}f(\tau, \o)d\xa(\tau, \o)$ exists for all $f \in \C^{\beta}(\R)$ with $\beta > \frac{1-\alpha}{2}$. The paper \cite{Ln} contains a special case of this result for functions $f(\xa(\cdot, \o)) \in \C^{1}(\R)$ for a.e.\ $\o$, derived by expanding $f$ in Taylor series in $\xa$ and using the fact that the quadratic variation of $\xa$ is zero. 

Similarly,  Theorem ~\ref{sum} holds with probability one for $\int_{0}^{t}f(\tau, \o)d\xa(\tau, \o)$ for all $f \in \C^{\beta}(\R)$ with $\beta > \frac{1-\alpha}{2}$.

It\^o's formula for fBm can be stated pathwise, as a corollary of Lemma ~\ref{ito}, however it can be stated also under weaker assumptions.
It\^o's formula for fBm has been established under very different assumptions in 
\cite{DH} for $0 < \alpha <1,$ 
and in \cite{DU} for $-1 <\alpha< 1$. 
\begin{Lemma}\label{ito-fbm}
Let $u:[0,\infty)\times \R \rightarrow \R$, $e, f:[0,\infty)\times\O \rightarrow \R$,  $0 < \alpha <1$, $\beta > \frac{1-\alpha}{2}$, $T > 0$.
Suppose $f \in C^{\beta}([0,T])$, $e$ is continuous, $u$ is differentiable  with continuous $\partial u/\partial t$ and $\partial u/ \partial x \in C^{\gamma}([0,T])$.	
For $0 \leq t \leq T$, consider the stochastic integral
\be
	\eta(t,\o)=\eta(0,\o)+\int_{0}^{t}e(s,\o)ds+\int_{0}^{t}f(s,\o)d\xa(s,\o).
\ee
Suppose $\sup\limits_{0 \leq t \leq T} E\left(\frac{\partial^2 u}{ (\partial \eta)^2} (t, \eta(t,\o))\right)^2 < \infty$.
Then for all $0 \leq t \leq T$, $v(t)=u(t,\eta(t,\o))$ is also a stochastic integral whose differential is
\be
	dv = \frac{\partial u}{\partial t} dt + 
			\frac{\partial u}{\partial \eta} d\eta.
\ee
\end{Lemma}
Note that for $0 < \alpha <1$,  It\^o's formula for fBm is the same as in the deterministic case. This follows from the fact that the quadratic variation of fBm is zero.
Theorem ~\ref{main} implies the following existence-uniqueness theorem for SDE's with fBm:
\begin{Theorem}\label{fbm_main}
Let $b, \sigma : [0, \infty) \times \R \rightarrow \R$,  $Z:{\bf \O} \rightarrow \R$ and $0 < \alpha < 1$, $\frac{1-\alpha}{2} < \beta < \frac{1+\alpha}{2}$. 
Suppose $b$ is globally Lipshitz in $t$ and $x$,
and $\sigma \in \C^{1}(\R)$ with $\sigma$, $\sigma'_t$, $\sigma'_x$ globally Lipshitz in $t$ and $x$.
Then for every $T > 0$ the SDE
\begin{equation}\label{fbm_sde}
	dX(t, \o) = b(t, X(t, \o))dt
	+ \sigma(t, X(t, \o))d\xa(t, \o), \quad X(0, \o) = Z(\o)  
\end{equation}
has a unique solution in $C^{\beta}([0,T])$ with probability 1.  
\end{Theorem}

In \cite{Ln, DH}, the existence and uniqueness theorem for 
solutions of stochastic differential equations  was proved 
when the diffusion coefficient is a function of time $t$ only;
in \cite{Ln}, the existence theorem was proved also when
the drift and diffusion coefficients are functions of $X$ only. Both papers adapt the methods used for oBm. The new idea in this paper is to use the formula in Theorem~\ref{sum}, which allows us to prove the existence and uniqueness theorem in the general case, when the drift and diffusion are functions of both $t$ and $X$.
\section{Additional results for fractional Brownian motion.}\label{additional}

When the integrand is of the form $f(\tau,\xa(\tau, \o))$, we can obtain 
estimates of the tail of the maximum of stochastic 
integrals from Theorem~\ref{sum} and from tail estimates  of the H\"older coefficient of fBm.
\begin{Theorem}\label{max-of-f}
Let $f: [0, 1]\times \R \rightarrow \R$  and $\frac{1}{2} < \gamma < \frac{1+\alpha}{2}$, $\delta = \frac{1+\alpha}{2} - \gamma$. Suppose 
$f$ is differentiable with bounded $|f'_t(t,x)|$, $|f'_x(t,x)|$. Write $|f'_t| = \sup\limits_{[0,1],\R}|f'_t(t,x)|$ and $|f'_x| = \sup\limits_{[0,1],\R}|f'_x(t,x)|$. Then
\bee\label{max-f}
\qquad	P\bigl(\max\limits_{0\leq t\leq 1} 
	\int_{0}^{t} f(\tau, \xa(\tau,\o)) d\xa(\tau, \o)
	> \lambda \bigr) \leq 
	\frac{2^{\gamma}+1}{2^{\gamma}-1} \sqrt{\frac{2}{\pi}}
		\sum_{n=1}^{\infty} \frac{2^{(1-\delta)n}}{\nu} 
		\exp\{-\frac{2^{\gamma}-1}{2^{\gamma}+1}
				\nu^2 2^{2n\delta-1}\}.
\eee
where 
\bee\label{nu}
	\nu = \frac{2^{2\gamma-1}-1}{|f'_x|}\left[
	\sqrt{(	|f(0,0)| + \frac{|f'_t|}{2^{\gamma+1}-2})^2
		+ \frac{4}{2^{2\gamma}-2}|f'_x|\lambda}
		- |f(0,0)| - \frac{|f'_t|}{2^{\gamma+1}-2}\right].
\eee
\end{Theorem}
	More information on the stochastic integral is available
when $f$ is a function of $t$ only, since in this case the techniques used 
for oBm can be applied.
\begin{Statement}\label{Aha}

Let $0 < \alpha < 1$, and let
$f(t, \o) = f(t)$ be a function of $t$ only. 
If $f \in L^{\frac{2}{1+\alpha}}([0, \infty))$, then 
the fBm stochastic integral $\int_{0}^{t}f(\tau)d\xa(\tau, \o)$ 
exists in $L^2([0, \infty) \times \O)$ for all 
$t \in [0, \infty)$.  

\end{Statement}
{\noindent}The proof is based on the Hardy-Littlewood-Sobolev inequality (see \cite{LL}).

Since $\int_{0}^{t}f(\tau)d\xa(\tau, \o)$ is a Gaussian process, we can show the following properties:
\begin{Statement}\label{max}

Let $\alpha$ and $f$ be as in Theorem~\ref{Aha}, and let
$0 < \beta < \alpha$.  Write 
$q_f(s, t) 
= \int_{s}^{t}\int_{s}^{t} f(u)f(v)\frac{dudv}{|u-v|^{1-\alpha}}$.
If $f \in L^{\frac{2}{1+\beta}}([0,1])$, then
\begin{enumerate}

\item for a.e.\ $\o$, $\int_{0}^{t} f(\tau)d\xa(\tau, \o)$ 
has a $t$-continuous version for all $t\in [0, 1]\,${\rm ;}

\item for every real $r$, 
$P\bigl(\max\limits_{0\leq t\leq 1} \int_{0}^{t} f(\tau) d\xa(\tau, \o)
> \lambda \bigr)$ is bounded from above by
$$
\int_{\lambda r/{\sqrt{q_{f_+}(0, 1)}}}^{\infty}  + 
\int_{\lambda (1-r)/{\sqrt{q_{f_-}(0, 1)}}}^{\infty}
\sqrt{\frac{2}{\pi}}\, e^{-x^2/2} dx, \quad 
{\rm~where~}
\quad  f_{\pm} = \frac{|f|\pm f}{2}\, ;
$$

\item for every integer $m\geq 2$ and real 
$\lambda\geq \sqrt{1 + \log{m^4}}$,
$P\bigl(\max\limits_{0\leq t\leq 1} \left|\int_{0}^{t} f(\tau)
d\xa(\tau, \o)\right| > \lambda \bigr)$ is bounded from above by
$$
\int_{\lambda/c}^{\infty} \frac{5}{2}\, m^2\, e^{-x^2/2} dx,
\quad
{\rm~where~}
\quad
c = \sup\limits_{0\leq s, t\leq 1} \sqrt{q_f(s, t)} + 
(2 + \sqrt{2})\int_1^{\infty} \sup\limits_{|s-t| < m^{-x^2}} \sqrt{q_f(s, t)}dx
< \infty.
$$

\end{enumerate}

\end{Statement}
{\noindent}(1) follows from Kolmogorov's continuity criterion, the bound (2) follows from Slepian's lemma \cite{Sl} (the interesting case is $0 < r < 1$) , and (3) follows from Fernique's inequality \cite{F}.

\subsection{Acknowledgments.}
I would like to thank my advisor Yakov Sinai for suggesting the problem and for all his help during the course of the work, Almut Burchard for suggesting a way to simplify the argument of section~\ref{sec:local}, the referee for his helpful suggestions in restructuring the paper and especially for an elegant way to simplify the argument of section~\ref{sec:extension}, Tadashi Tokieda for his help and comments on the text and Luc Rey-Bellet for his comments on the text.

\section{Proofs of theorems.}
\subsection{Proof of Theorem~\ref{sum}.}

\noindent 
{\em Proof:} 
	The idea of the proof is to write the Riemann sums on the smaller scales in terms of the Riemann sums on the larger scales as in Renormalization Group.
	Denote by $S^n(f)$ the Riemann sum of $f$ corresponding to 
the partition of $[0, 1]$ into $2^n$ equal sub-intervals.  We
have
\be
S^n(f)    
   &=&	S^{n-1}(f) + 
	\sum_{i=0}^{2^{n-1}-1}
	\Delta f(s_{2i}^n) \Delta g(s_{2i+1}^n) \\
   &=&   \cdots 							  \\
   &=&   S^0(f) + 
	\sum_{k=1}^n \sum_{i=0}^{2^{k-1}-1}
	\Delta f(s_{2i}^k) \Delta g(s_{2i+1}^k).
\ee
As $n \to \infty$, $S^n(f)$ converges to 
$\int_{0}^{t} f(\tau)dg(\tau)$ by Theorem~\ref{Stieltjes}. The right hand side converges provided $\beta+\gamma>1$. \qed

\subsection{Proof of the change of variables formula.}
The proof given here follows \cite{Mk} for the most part.

\noindent 
{\em Proof:} We can write the integral version of equation ~(\ref{ito1}):
\begin{displaymath}
	v(t)-v(0)= \int_{0}^{t}\frac{\partial u}{\partial s} ds +
		 	\int_{0}^{t} \frac{\partial u}{\partial \eta} d\eta
\end{displaymath}
From Theorem ~\ref{sum} and from continuity of $e$ it follows that:
\begin{displaymath}
	\Delta \eta(t) = f(t)\Delta g(t) +
			 e(t)\Delta t +
			 O(\Delta t)^{\beta+\gamma} + o(\Delta t).
\end{displaymath}
\begin{eqnarray*}
&&	v(t,g(t))-v(0) = \sum_{k \leq [2^n t]} \{u(\frac{k}{2^n},
	\eta (\frac{k}{2^n})) -
	u(\frac{k-1}{2^n}, 
		\eta (\frac{k}{2^n}))\} 		\cr 
 && \qquad +  	 \sum_{k \leq [2^n t]} \{u(\frac{k-1}{2^n},
	\eta (\frac{k}{2^n})) -
	u(\frac{k-1}{2^n}, \eta (\frac{k-1}{2^n}))
	\} + u(t,\eta(t))-
	u(\frac{[2^nt]}{2^n}, \eta (\frac{[2^nt]}{2^n})) \cr
 && \qquad =  	\sum_{k \leq [2^nt]}\{\frac{\partial u}{\partial t}
	(\frac{k-1}{2^n},
	\eta (\frac{k}{2^n}))\frac{1}{2^n}
		+o(\frac{t}{2^{n}})\}				\cr  
 && \qquad + 	\sum_{k \leq [2^nt]}\{\frac{\partial u}{\partial \eta}
	[\frac{k-1}{2^n},\eta (\frac{k-1}{2^n})]
	(\eta (\frac{k}{2^n})
		-\eta (\frac{k-1}{2^n}))   \cr  
 && \qquad + 	\frac{1}{2}\frac{\partial^2 u}{\partial \eta^2}
	(\frac{k-1}{2^n},\eta (\frac{k-1}{2^n}))
	(\eta (\frac{k}{2^n})
		-\eta (\frac{k-1}{2^n}))^2 + 
	o(\Delta \eta^2 ) \} 		 \cr   
 && \qquad = 	\sum_{k \leq [2^nt]}\{(\frac{\partial u}{\partial t}
	(\frac{k-1}{2^n},
	\eta (\frac{k}{2^n}))				
   + 	\frac{\partial u}{\partial \eta}
	(\frac{k-1}{2^n},\eta (\frac{k-1}{2^n}))
	e[\frac{k-1}{2^n}])\frac{1}{2^n}+o(\frac{t}{2^{n}})\} \cr
 && \qquad + 	\sum_{k \leq [2^nt]}\{\frac{\partial u}{\partial \eta}
	(\frac{k-1}{2^n},\eta (\frac{k-1}{2^n}))
	f(\frac{k-1}{2^n})
	( g(\frac{k}{2^n})
		-g(\frac{k-1}{2^n}))  \cr
 && \qquad + 	 \frac{1}{2}\frac{\partial^2 u}{\partial \eta^2}
	(\frac{k-1}{2^n},\eta (\frac{k-1}{2^n}))
	f(\frac{k-1}{2^n})^2
	(g(\frac{k}{2^n})
		-g(\frac{k-1}{2^n}))^2 + 
	o(\Delta g^2 )\} 				\cr
&& \qquad=	\int_{0}^{t}\left(\frac{\partial u}{\partial s}+
	\frac{\partial u}{\partial \eta}e(s)\right) ds + o(1) +
	\int_{0}^{t}\frac{\partial u}{\partial \eta} f(s) dg (s) +
	\frac{1}{2}\sum_{k \leq [2^nt]}u_{\eta\eta}f^{2} \Delta g^2 + 
	o\left(\sum \Delta g^2 \right).
\end{eqnarray*}
Since $ u_{\eta\eta}f^{2}\Delta g^2 \leq O\left(\frac{t}{2^{n}}\right)^{2\gamma},$
and $\gamma>\frac{1}{2}$, $\sum_{k \leq [2^nt]}u_{\eta\eta}f^{2}\Delta g^2 \rightarrow 0$ as $n \rightarrow \infty$. 
Similarly, $o\left(\sum \Delta g^2 \right) \rightarrow 0$ as $n \rightarrow \infty$. \qed 

\subsection{Proof of Theorem~\ref{main}: local existence and uniqueness.} 
\label{sec:local}

In this section we will prove the local existence and uniqueness result. We will derive the global result in Theorem~\ref{main} by showing that we can apply the local existence and uniqueness result repeatedly, taking as initial condition the value of the solution at the end of the previous interval. 

	Fix $T > 0$, $s \in [0, T]$ and $a \in \R $.  Define an integral operator
$$
	F_{t} X = \int_s^t b(\tau, X(\tau)) d\tau + 
		      \int_s^t\sigma(\tau, X(\tau)) dg(\tau) + a.
$$
For $s=0$ and $a = x_{0}$, solutions of the ODE~(\ref{sde}) are exactly fixed points of $F$.

	For $0 < \beta < 1$, and for $\ve > 0$, consider
the Banach space $\C^{\beta}([s,s+ \ve])$ with the norm $||f||_{\beta} = \max\limits_{t}|f(t)| + \max\limits_{t_1 \neq t_2}
			\frac {|f(t_1) - f(t_2)|}{|t_1 -t_2|^{\beta}}$. For $K > 0$, consider the closed subset 
$\C_K^{\beta}([s,s+ \ve])
= \{ f : |f(t_2) - f(t_1)| \leq K |t_2-t_1|^{\beta},\  \forall t_1, t_2 \in [s, s+\ve] \}$ of $\C^{\beta}([s,s+ \ve])$.
Finally, consider the closed subset $\C_K^{\beta}([s,s+ \ve], a) = \{ f \in \C_K^{\beta}([s,s+ \ve]) : f(s) = a \}$ of $\C^{\beta}([s,s+ \ve])$. 
For a given $g \in \C^{\gamma}([0,T])$, there is $L > 0$ such that $g \in \C_{L}^{\gamma}([0,T])$.
	We will show by the contraction mapping theorem 
that for given $T$, $s$, $a$, $K$, $L$,  
there exists $\ve > 0$ 
such that $F_t$ has a unique fixed point on $\C_K^{\beta}([s,s+\ve], a)$.  
We need to establish that $F_t$ maps 
$\C_K^{\beta}([s,s+\ve], a)$ into itself and that it is a contraction. 
This will be done in Lemma~\ref{mapitself}. The necessary preliminary estimates are obtained in Corollaries~\ref{corollary} and~\ref{corollary_of_corollary} of Theorem~\ref{sum}. In what follows we will consider $T, K, L > 0$ to be fixed.
\begin{Corollary}\label{corollary}
Let $\beta$ and  $\gamma$  be as in Theorem~\ref{Stieltjes}, $T>0$ and let $g \in \C_L^{\gamma}([0,T])$. Then for every $s,t \in [0,T]$,
\begin{equation}\label{corollary_sum}
||\int_s^t X(\tau)dg(\tau)||_{\infty} 
	= L ||X||_{\beta} (t-s)^{\gamma}(1 		
	+
	 \frac{(t-s)^{\beta}}{2^{\beta+\gamma}-2}).
\end{equation}
\end{Corollary}

\noindent 
{\em Proof:} 
By Theorem ~\ref{sum} 
\begin{equation}\label{corollary_sum1}
|\int_s^t X(\tau)dg(\tau)| \leq |X (s)| |g (t) - g (s)| + \sum_{k=1}^{\infty}\sum_{i=0}^{2^{k-1}-1}|\Delta X(s + s_{2i}^{k})| |\Delta g (s + s_{2i+1}^{k})|.
\end{equation}
Since $|X (s)| \leq ||X||_{\beta}$, $|\Delta X(s + s_{2i}^{k})| \leq||X||_{\beta}(\frac{t-s}{2^k})^{\beta}$ and $g \in \C_L^{\gamma}([0,T])$, we obtain after performing the sum in~(\ref{corollary_sum1}),
\begin{displaymath}
|\int_s^t X(\tau)dg(\tau)| \leq L ||X||_{\beta} (t-s)^{\gamma} + 
		L ||X||_{\beta}\frac{(t-s)^{\beta+\gamma}}{2^{\beta+\gamma}-2}.
\end{displaymath}
\qed

\begin{Corollary}\label{corollary_of_corollary}
Let $\ve  > 0$ and let $\beta$, $\gamma$, $g$ and $T$ be as in Corollary~\ref{corollary}. Then, for every $s \in [0,T]$ and $t \in [s, s+\ve]$,
\begin{equation}\label{corollary_of_corollary1}
||\int_s^t X(\tau)d\tau||_{\beta} \leq ||X||_{\infty} \ve^{1-\beta}(1+\ve^{\beta}),
\end{equation}
and
\begin{equation}\label{corollary_of_corollary2}
||\int_s^t X(\tau)dg (\tau)||_{\beta} \leq L ||X||_{\beta} \ve^{\gamma-\beta}(1+\ve^{\beta})(1+\frac{\ve^{\beta}}{2^{\beta+\gamma}-2}),
\end{equation}
\end{Corollary}

\noindent 
{\em Proof:} 
\begin{eqnarray*}
||\int_s^t X(\tau)d\tau||_{\beta} &=& \max\limits_{t \in [s, s + \ve]}|\int_s^t X(\tau)d\tau| + \max\limits_{t_{1} \not= t_{2} \in [s, s + \ve]}\frac{|\int_{t_{1}}^{t_{2}} X(\tau)d\tau|}{|t_{2}-t_{1}|^{\beta}} \\
&\leq& ||X||_{\infty} (\ve + \ve^{1-\beta}).
\end{eqnarray*}
From Corollary~\ref{corollary} we obtain
\begin{equation}\label{corollary_of_corollary3}
\max\limits_{t \in [s, s + \ve]}|\int_s^t X(\tau)dg (\tau)| \leq L ||X||_{\beta}\ve^{\gamma}(1+\frac{\ve^{\beta}}{2^{\beta+\gamma}-2}),
\end{equation}
similarly we obtain
\begin{equation}\label{corollary_of_corollary4}
\max\limits_{t_{1} \not= t_{2} \in [s, s + \ve]}\frac{|\int_{t_{1}}^{t_{2}} X(\tau)dg (\tau)|}{|t_{1}-t_{2}|^{\beta}} \leq L ||X||_{\beta}\ve^{\gamma-\beta}(1+\frac{\ve^{\beta}}{2^{\beta+\gamma}-2}).
\end{equation}
The result~(\ref{corollary_of_corollary2}) follows from ~(\ref{corollary_of_corollary3}) and ~(\ref{corollary_of_corollary4}).
\qed
\begin{Lemma}\label{mapitself}
Let $b, \sigma : [0, \infty) \times \R \rightarrow \R$, $1 \geq \gamma > 1/2$, $T,\ K,\ L > 0$. 
Suppose $b$ and $\sigma$ are  Lipschitz in $x$ and $t$.
Then there exists $\ve _{1} > 0$, such that for all $
\gamma > \beta > 1-\gamma$ and $t \in [s, s+\ve_1],$ the operator 
$F_t$ maps $\C_K^{\beta}([s,s+\ve_1], a)$ into itself.
\end{Lemma}
\noindent {\it Notation.} 
We shall denote the Lipschitz coefficients of $b$ and $\sigma$ by ${\bf B}$ and ${\bf S}$ respectively. 

\noindent 
{\em Proof:}  
Let $\ve > 0$ and $t_1, t_2 \in [s, s+\ve]$. It is sufficient to demonstrate that $||F_{t}X||_{\beta}\leq K$ for a sufficiently small $\ve >0$.

By  the triangle inequality and by Corollary ~\ref{corollary_of_corollary},
\begin{equation}\label{mapitself1}
	||F_{t}X||_{\beta} \leq \ve^{\gamma-\beta} (1+\ve^{\beta})(||b(t,X(t))||_{\infty} \ve^{1-\gamma}+ L ||\sigma(t,X(t))||_{\beta} (1+\frac{\ve^{\beta}}{2^{\beta+\gamma}-2})).
\end{equation}
Since $b$ and $\sigma$ are Lipschitz and $X \in \C_K^{\beta}([s,s+\ve], a)$, it is easy to see that
\begin{equation}\label{mapitself2}
	||b(t,X(t))||_{\infty} \leq |b (s,a)| + {\bf B}(\ve + K\ve ^{\beta}),
\end{equation}
and
\begin{equation}\label{mapitself3}
	||\sigma(t,X(t))||_{\beta} \leq |\sigma(s,a)| + {\bf S} (\ve ^{1-\beta}+ K) (1 + \ve ^{\beta}).
\end{equation}
Substituting~(\ref{mapitself2}) and~(\ref{mapitself3}) into~(\ref{mapitself1}), we obtain
\begin{eqnarray}\label{mapitself4}
	||F_{t}X||_{\beta} &\leq& \ve ^{\gamma-\beta}(1+\ve^{\beta})\big
[|b (s,a)| + {\bf B}(\ve + K\ve ^{\beta}))\ve ^{1-\gamma} \cr && + L(1+\frac{\ve^{\beta}}{2^{\beta+\gamma}-2}) (|\sigma(s,a)| + {\bf S} (1 + \ve ^{\beta})(\ve ^{1-\beta}+ K) )\big].
\end{eqnarray}

Since the right hand side is an increasing continuous function 
of $\ve$ and is $0$ at $\ve = 0$, it equals $K$ at some $\ve_1$.  For this choice of $\ve_1$
(or any smaller $\ve_1$), $F_t$ maps $\C_K^{\beta}([s, s+\ve_1], a)$ into itself. \qed

\begin{Lemma}\label{contraction}
Assume the same hypothesis as in Lemma~\ref{mapitself}. 
Suppose $b$ is Lipschitz in $t$ and $x$ and $\sigma \in \C^{1} ([0,\infty) \times \R)$ with  
$\sigma'_t(t,x)$, $\sigma'_x(t,x)$ Lipschitz in $x$.
Then there exists $\ve_2 > 0$ such that for all $
\gamma > \beta > 1-\gamma$ and $t \in [s,s+\ve_2]$, 
the operator $F_t$ is a contraction on 
$\C_K^{\beta}([s,s+\ve_2])$. 
\end{Lemma}
\noindent {\it Notation.} 
We shall denote the Lipschitz coefficient of $b$ by ${\bf B}$ and the Lipschitz coefficient of $\sigma$, $\sigma'_t$, $\sigma'_x$ by ${\bf S}$.

\noindent 
{\em Proof:}  
We need to show that there exist $\ve_2 > 0$ and $\lambda < 1$ such
that for all $t \in [s,s+\ve_2]$ and all $X, Y \in \C_K^{\beta}([s,s+\ve_2])$,
\be
	||F_tX- F_tY||_{\beta} \leq \lambda ||X-Y||_{\beta}. 
\ee
By  the triangle inequality and by Corollary ~\ref{corollary_of_corollary},
\begin{eqnarray}\label{contraction1}
	||F_{t}X- F_tY||_{\beta} &\leq& \ve^{\gamma-\beta} (1+\ve ^{\beta})\big[||b(t,X(t))-b (t,Y (t)) ||_{\infty} \ve ^{1-\gamma}\cr &&+ L ||\sigma(t,X(t))-\sigma(t,Y(t))||_{\beta} (1+\frac{\ve^{\beta}}{2^{\beta+\gamma}-2})\big].
\end{eqnarray}
We estimate the two terms in~(\ref{contraction1}) separately.
Since $b$ is Lipschitz,
\begin{equation}\label{contraction2}
||b(t,X(t))-b (t,Y (t)) ||_{\infty} \leq {\bf B}|X (t)-Y (t)| \leq {\bf B}||X-Y||_{\beta}.
\end{equation}
Now we will estimate $||\sigma(t,X(t))-\sigma(t,Y(t))||_{\beta}$. Since $\sigma$ is differentiable,
\begin{equation}\label{contraction3}
	\max_{t} |\sigma(t,X(t))-\sigma(t,Y(t))| \leq {\bf S} ||X-Y||_{\beta}.
\end{equation}
By the fundamental theorem of calculus,
\begin{displaymath}
	\sigma(t,X(t))-\sigma(t,Y(t)) = (X (t) - Y (t))\int_{0}^{1}\sigma'_x(t, \nu X(t)+(1-\nu) Y(t)) d\nu.
\end{displaymath}	
Therefore
\be
&& |\sigma(t_{1},X(t_{1}))-\sigma(t_{1},Y(t_{1}))-\sigma(t_{2},X(t_{2}))+\sigma(t_{2},Y(t_{2}))| = \cr &&
	|(X(t_{1}) - Y(t_{1})
	- X(t_{2}) + Y(t_{2}))
	\int_{0}^{1}\sigma'_x(t_{1}, \nu X(t_{1})+
	(1-\nu) Y(t_{1})) d\nu \cr &&
	+ (X(t_{2}) - Y(t_{2})) 
	\int_{0}^{1}\left(\sigma'_x(t_{1}, \nu X(t_{1})+
	(1-\nu) Y(t_{1})) -
	\sigma'_x(t_{2},\nu X(t_{2})+
	(1-\nu) Y(t_{2}))\right) d\nu| \cr &&
	\leq ||X-Y||_{\beta} (t_{2}-t_{1})^{\beta}{\bf S} + ||X-Y||_{\beta}
	\left({\bf S} (t_{2}-t_{1}) + {\bf S} (							\nu |X(t_{1})-X(t_{2})| + (1-\nu)|Y(t_{1})-Y(t_{2})|)\right)  \cr &&
	\leq ||X-Y||_{\beta} (t_{2}-t_{1})^{\beta}{\bf S} (1+(t_{2}-t_{1})^{1-\beta} + K).
\ee
Consequently
\begin{equation}\label{contraction4}
	||\sigma(t,X(t))-\sigma(t,Y(t))||_{\beta} \leq ||X-Y||_{\beta} {\bf S} (2+\ve ^{1-\beta} + K).
\end{equation}
Substituting~(\ref{contraction3}) and~(\ref{contraction4}) in~(\ref{contraction1}), we obtain 
\begin{equation}\label{contraction5}
	||F_{t}X- F_tY||_{\beta}\leq \ve^{\gamma-\beta} (1+\ve ^{\beta})\big[{\bf B} \ve ^{1-\gamma}+ L {\bf S}(2+\ve ^{1-\beta} + K) (1+\frac{\ve^{\beta}}{2^{\beta+\gamma}-2})\big] ||X- Y||_{\beta}.
\end{equation}
The coefficient of $||X- Y||_{\beta}$ is an increasing function of 
$\ve$ and is $0$ at $\ve = 0$.  Choose $\ve_2 > 0$ small enough so that this coefficient is 
less than 1 and $\ve_2 \leq \ve_1$.  Then $F$ is a contraction on 
$\C_K^{\beta}([s,s+\ve_2]).$
\qed

Combining Lemmas~\ref{mapitself} and~\ref{contraction}, we obtain a local existence and uniqueness result:
%
%

%
%
\begin{Corollary}\label{local-existence-uniqueness}
Assume the same hypothesis as in Lemma~\ref{contraction}.
{\noindent}Then there exists $\ve > 0$ depending on $s$ such that for all $t \in [s, s+\ve]$, ODE~(\ref{sde}) has a unique solution in $\C_K^{\beta}([s,s+\ve], a)$. 
In particular, if ODE ~(\ref{sde}) has a solution $X$ on  $[0,s]$, then there exists $\ve > 0$ depending on $s$ such that for all $t \in [s, s+\ve]$, ODE~(\ref{sde}) has a unique solution $Y$ in $\C_K^{\beta}([s,s+\ve], X(s))$.
\end{Corollary}

\noindent
{\em Proof:}  
 Take $\ve$ to be $\ve_2$ of Lemma~\ref{contraction}.
Since $F_t$ is a contraction on the closed subset $\C_K^{\beta}([s,s+\ve], a)$ of the complete metric space 
$\C^{\beta}([s,s+\ve])$, it has a unique fixed point $X$ in
$\C_K^{\beta}([s,s+\ve], a)$.  From the definition of $F_t$ 
it follows that $X$ is a unique solution of ODE~(\ref{sde}) on $[s,s+\ve]$ in $\C_K^{\beta}([s,s+\ve], a)$. 
\qed

The sufficient conditions on $\ve$ in Corollary~\ref{local-existence-uniqueness} are given by inequalities~(\ref{mapitself4}) and~(\ref{contraction5}) with $a$ replaced by $X(s)$:
\begin{eqnarray}\label{sufficient-condition1}
	&&\ve ^{\gamma-\beta}(1+\ve^{\beta})\big[|b (s,X (s))| + {\bf B}(\ve + K\ve ^{\beta}))\ve ^{1-\gamma}\cr && \qquad + L (1+\frac{\ve^{\beta}}{2^{\beta+\gamma}-2})(|\sigma(s,X (s))| + {\bf S} (1 + \ve ^{\beta})(\ve ^{1-\beta}+ K) )\big] \leq K.
\end{eqnarray}
\begin{equation}\label{sufficient-condition2}
	 \ve^{\gamma-\beta} (1+\ve ^{\beta})({\bf B} \ve ^{1-\gamma}+ L {\bf S}(2+\ve ^{1-\beta} + K) (1+\frac{\ve^{\beta}}{2^{\beta+\gamma}-2})) < 1.
\end{equation}
Inequality ~(\ref{sufficient-condition2}) does not depend on $X(s)$.  

\subsection{Proof of Theorem ~\ref{main}: global existence and uniqueness.}\label{sec:extension}

By Corollary ~\ref{local-existence-uniqueness} with $s=0$ and $a = x_{0}$, ODE  ~(\ref{sde}) has a unique solution on $[0, \ve_0]$, where $\ve_0$ satisfies~(\ref{sufficient-condition1})  and~(\ref{sufficient-condition2}) with $s=0$. 
By using Corollary~\ref{local-existence-uniqueness} $n$ times we obtain that the solution exists on $[0, \ve_0+\cdots+\ve_{n-1}]$. ODE~(\ref{sde}) has  a solution on $[0, T]$ if  there exists $m > 0$ such that $\sum_{i=0}^{m} \ve_i \geq T$. This is true if $b$ and $\sigma$ are globally bounded, since in this case we can choose $\ve _{i} = \ve _{0}$ (to see this we can substitute the global bounds on $b$, $\sigma$, ${\bf B}$, ${\bf S}$ into~(\ref{sufficient-condition1}) and~(\ref{sufficient-condition2})). In the case when $b$ and $\sigma$ grow at most linearly, we will use a change of variables to reduce it to the case of globally bounded $b$ and $\sigma$. 

\noindent 
{\em Proof:} 
{\it Existence\/}: 
Suppose that $b$ and $\sigma$ are bounded on $[0,\infty)\times \R $, then taking the upper bound on $b$ and $\sigma$ in~(\ref{sufficient-condition1}) we get that $\ve_i$ satisfying~(\ref{sufficient-condition1}) and~(\ref{sufficient-condition2}) does not depend on $i$. In this case the global existence is established.
Now suppose that $b$ and $\sigma$ satisfy the assumptions of Theorem~\ref{main}. Consider the ODE 
\begin{equation}\label{ode}
	dy(t) = \frac{b(t, \tan y (t))}{1+ ( \tan y (t))^{2}}dt
	+ \frac{\sigma(t,\tan y (t))}{1+ (\tan y (t))^{2}}d g(t). 
\end{equation}
This ODE has globally bounded coefficients satisfying the assumptions of Theorem~\ref{main}, and thus has a global solution on $[0,T]$. Now, $x (t)= \tan y (t)$ satisfies equation~(\ref{sde}) (by Lemma~\ref{ito}). Thus equation~(\ref{sde}) has a global solution on $[0,T]$. 

{\it Uniqueness\/}:
Let $Y_1$ and $Y_2$ be two solutions in 
$\C^{\beta}([0,T])$.  Then there exist $K_1$ and $K_2$ 
such that $Y_1 \in \C_{K_1}^{\beta}([0,T])$ and 
$Y_2 \in \C_{K_2}^{\beta}([0,T])$, so 
$Y_1$ and $Y_2$ are in $\C_{\max\{K_1, K_2\}}^{\beta}([0,T])$.
$Y_1$ and $Y_2$ coincide at the initial point $t = 0$. 
Let $t_{\sup}$ be the supremum of the 
set on which they coincide.  Since both solutions are continuous, 
they coincide at $t_{\sup}$ as well.  $t_{\sup}$ must 
equal $T$, for otherwise we can make $Y_1$ and $Y_2$ coincide past 
$t_{\sup}$ by Corollary~\ref{local-existence-uniqueness}.
\qed

\subsection{Proof of It\^o's formula for fBm (Lemma~\ref{ito-fbm}).}
\label{sec:Itoforfbm}

\noindent 
{\em Proof:} 
By analogy with the proof of Lemma ~\ref{ito}, we obtain
\be
       v(t,\xi(t))-v(0)=\int_{0}^{t}\left(\frac{\partial u}{\partial s}+
	\frac{\partial u}{\partial \eta}e(s)\right) ds + o(1) +
	\int_{0}^{t}\frac{\partial u}{\partial \eta} f(s) d\xa +
	\frac{1}{2}\sum_{k \leq [2^nt]}u_{\eta\eta}f^{2}\Delta\xa^2 + 
	o\left(\sum \Delta \xa^2 \right).
\ee
Since $ E(u_{\eta\eta}f^{2}\Delta\xa^2) \leq ||f||_{\infty}\sqrt{ E(u_{\eta\eta})^2}
	\sqrt{ E(\Delta\xa^4)} = O\left(\frac{t}{2^{n}}\right)^{2\gamma},$
by Chebyshev inequality
$$
	P(\sum_{k \leq [2^nt]} u_{\eta\eta}f^{2}\Delta\xa^2 \geq \frac{1}{n})
	\leq {\rm const~} \frac{n}{2^{n(2\gamma-1)}},
$$
and by the Borel-Cantelli lemma $\sum_{k \leq [2^nt]}u_{\eta\eta}f^{2}\Delta\xa^2 \rightarrow 0$ as $n \rightarrow \infty$ a.e. 
An analogous argument shows that $o\left(\sum \Delta \xa^2 \right) \rightarrow 0$ as $n \rightarrow \infty$ a.e.
Thus It\^o's formula holds.
\qed 
\subsection{Proof of Theorem~\ref{max-of-f}}\label{proofmax}
\noindent 
{\em Proof:} 
Since $f$ is differentiable, 
and since $\xa(\cdot,\o) \in C^{\gamma}([0,1])$ with probability 1, for a.e. $\o$ there exists $L(\o) > 0$ such that
$$|\Delta f(s_{2 i}^{k}, \xa(s_{2 i}^{k}))| \leq |f'_t| \frac{t}{2^k} + |f'_x| \frac{L(\o) t^{\gamma}}{2^{k \gamma}}$$ holds.
From Theorem ~\ref{sum} we get 
\be
	\left|\int_{0}^{t} f(\tau, \xa(\tau))d\xa(\tau)\right| \leq
	|f(0,0)| L(\o)t^{\gamma}+ |f'_t|\frac{ L(\o)t^{\gamma+1}}{2^{\gamma+1}-2} 
		+ |f'_x|\frac{ L(\o)^2t^{2 \gamma}}{2^{2\gamma}-2},
\ee
and so,
\be
	\max\limits_{0\leq t\leq 1} 
	\left|\int_{0}^{t} f(\tau, \xa(\tau))d\xa(\tau)\right| \leq
	|f(0,0)| L(\o) + |f'_t|\frac{ L(\o)}{2^{\gamma+1}-2} + 
		|f'_x|\frac{ L(\o)^2}{2^{2\gamma}-2}.
\ee
Therefore 
\be
&&	P\{\o: \max\limits_{0\leq t\leq 1} 
	\left|\int_{0}^{t} f(\tau, \xa(\tau))d\xa(\tau)\right| > \lambda \} \\ 
	&&\leq P\{\o: |f(0,0)| L(\o) + |f'_t|\frac{ L(\o)}{2^{\gamma+1}-2} + 
		|f'_x|\frac{ L(\o)^2}{2^{2\gamma}-2} > \lambda \} \\
	&&\leq P\{\o: L(\o)  > \nu \} \\
	&&\leq P\{\o: \exists\ t_1, t_2 \in [0,1] \ {\rm s.t.} \ 
			|\xa(t_2) - \xa(t_1)| > \nu |t_2 - t_1|^{\gamma}\},
\ee
where $\nu$ is given by ~(\ref{nu}).

It is  easy to see that if
\be
		|\xa(\frac{k+1}{2^n}) - \xa(\frac{k}{2^n})| \leq 
		\frac{2^{\gamma}-1}{2^{\gamma}+1} \frac{L(\o)}{2^{n \gamma}}, \quad				\forall\ n > 0, \ 0 \leq k \leq 2^n -1, 
\ee
then
\be
	 |\xa(t_2) - \xa(t_1)| \leq 	
		L(\o) |t_2 -t_1|^{\gamma}, \quad \forall\ t_1, t_2 \in [0,1].
\ee
Therefore
\be
&&	P\{\o: \exists\ t_1, t_2 \in [0,1] \ {\rm s.t.} \ 
			|\xa(t_2) - \xa(t_1)| > \nu |t_2 - t_1|^{\gamma}\} \\
	&\leq& P\{\o: \exists\ n > 0, 0 \leq k \leq 2^n -1, \ {\rm s.t.} \
		 |\xa(\frac{k+1}{2^n}) - \xa(\frac{k}{2^n})| >
		\frac{2^{\gamma}-1}{2^{\gamma}+1} \frac{\nu}{2^{n \gamma}}\} \\
	&\leq& \sum_{n=1}^{\infty} \sum_{k=0}^{2^n-1} 
		P\{\o:	|\xa(\frac{k+1}{2^n}) - \xa(\frac{k}{2^n})| > 
		\frac{2^{\gamma}-1}{2^{\gamma}+1} \frac{\nu}{2^{n \gamma}}\} \\
	&\leq& \sum_{n=1}^{\infty} 2^n \sqrt{\frac{2}{\pi}} 
		\int_{\frac{2^{\gamma}-1}{2^{\gamma}+1} \nu 2^{n \delta}}^{\infty} 
							e^{-x^2/2} dx.
\ee
Using the estimate $\int_{c}^{\infty} e^{-x^2/2} dx \leq \frac{e^{-c^2/2}}{c}$, we obtain the desired result.	
\qed
\subsection{Proof of Statement~\ref{Aha}.}

\noindent 
{\em Proof:} 

The proof will be reached via the step-by-step procedure used for oBm.
Fix $t > 0$.

\noindent		
{\it Step 1}. 
	Let $\phi: [0, \infty) \rightarrow \R$  be a simple function of the 
form $\sum_{i=0}^{2^n-1} \phi(s^n_i)\chi_{[s^n_i, s^n_{i+1}]}$, where 
$\chi$ is an indicator function and $s^n_i = \frac{i}{2^n}t$.
Define
$$
	\int_0^t \phi(\tau)d\xa(\tau) =
		\sum_{i=0}^{2^n-1} \phi({s}^n_{i})
	 \Delta \xa({s}^n_{i}).
$$

\noindent
{\it Step 2}. 
	Let $g \in \C^{\infty}([0, \infty))$.  Approximate $g$ by a 
sequence of simple functions: 
$\phi_n(\tau) = \sum_{i=0}^{2^n-1} g(s^n_i)\chi_{[s^n_i, s^n_{i+1}]}(\tau)$.
Then $\phi_n \rightarrow g$ uniformly on $[0, t]$ and 
$\int_0^t |g(\tau) - \phi_n(\tau)|^{\frac{2}{1+\alpha}} d\tau \rightarrow 0$. 
Therefore the sequence $\phi_n$ is Cauchy in $L^{\frac{2}{1+\alpha}}([0, \infty))$.
To show that the sequence $\int_0^t \phi_n(\tau)d\xa(\tau)$ is 
Cauchy in $L^{2}([0, \infty) \times \O)$, we use the 
Hardy-Littlewood-Sobolev inequality: 
\be
&&	E\left(\int_0^t \phi_m(\tau)d\xa(\tau) -
		\int_0^t \phi_n(\tau)d\xa(\tau)\right)^{\! 2} \\		
	&& \qquad = C \alpha(\alpha + 1) \int_0^t \int_0^t
	\frac {(\phi_m(u)-\phi_n(u)) (\phi_m(v)-\phi_n(v))}
		{|u-v|^{1-\alpha}} du dv \\
	&& \qquad \leq {\rm~const~}\cdot
		\|\phi_m - \phi_n\|^2_{\frac{2}{1+\alpha}}
	\rightarrow 0	\quad {\rm as} \qquad m,n \rightarrow \infty.
\ee
Thus the integral $\int_0^t g(\tau)d\xa(\tau)$ exists as the
$L^2$-limit of $\int_0^t \phi_n(\tau)d\xa(\tau)$.

\noindent
{\it Step 3}. 
	Let $f \in L^{\frac{2}{1+\alpha}}([0, \infty))$. 
Let $j \in \C^{\infty}_c([0, \infty))$ with $\int_{[0, \infty)}j = 1$. 
Define $j_{n}(\tau) = \frac{1}{n} j(n \tau)$ and 
$g_n = j_{n}*f$.  Then $g_n \in \C^{\infty}([0, \infty))$ and
$\int_0^t |g_n(\tau)-f(\tau)|^{\frac{2}{1+\alpha}} d\tau \rightarrow 0$.  
In particular, the sequence $g_n$ is Cauchy in
$L^{\frac{2}{1+\alpha}}([0, \infty))$. 
The  Hardy-Littlewood-Sobolev inequality gives
\be
	E\left(\int_0^t g_m(\tau)d\xa(\tau)-
		\int_0^t g_n(\tau)d\xa(\tau)\right)^{\! 2} 	
	&\leq&  {\rm~const~}\cdot
		\|g_m - g_n\|^2_{\frac{2}{1+\alpha}}
	\rightarrow 0	\qquad {\rm as} \quad m, n \rightarrow \infty.	
\ee
Thus the integral $\int_0^t f(\tau)d\xa(\tau)$ exists as an $L^2$-limit of $\int_0^t g_n(\tau)d\xa(\tau)$.

Thus, for all $f \in L^{\frac{2}{1+\alpha}}([0,1])$, 
we can choose simple functions $\phi_n$
converging in $L^{\frac{2}{1+\alpha}}$ to $f$ such that the $L^2$-limit of
$\int_0^t \phi_n(\tau)d\xa(\tau)$ exists. 
\qed

\subsection{Proof of Statement~\ref{max}.}

\subsection{Lemmas}

	We begin with three lemmas.  The first is Slepian's
lemma \cite{Sl, Kl}. 

\begin{Lemma}\label{Slepian} (\it Slepian)
	Let $\Gamma$ be a countable set, and let $X(t)$,
$Y(t)$ be two real Gaussian processes indexed by $t\in \Gamma$.
Suppose $E X^2(t) = E Y^2(t)$ and $E X(s)X(t) \geq E Y(s)Y(t)$
for all $s, t \in \Gamma$.  Then, for all real $\lambda$,
$P\bigl(\max\limits_{t \in \Gamma} X(t) \geq \lambda\bigr)
	\leq  	
	P\bigl(\max\limits_{t \in \Gamma} Y(t) \geq \lambda\bigr).$
\end{Lemma}

\noindent
Consequently, if $X$ and $Y$ have continuous versions, 
Lemma~\ref{Slepian} holds when the index set $\Gamma$ is $[0,1]$. 

	The next lemma gives us Markov property.

\begin{Lemma}\label{Ymarkov-continuous}

Let $\alpha$, $f$ be as in Theorem~\ref{Aha}, and let
$0 < \beta < \alpha$.
Let $Y(t)$ be a Gaussian process such that $E Y(t) = 0$ and 
$E Y(s)Y(t) = q_f(0, s)
	= \int_0^s\int_0^s f(u)f(v)\frac{dudv}{|u-v|^{1-\alpha}}$ 
does not depend on $t$ whenever $s \leq t$. 
If $f \in L^{\frac{2}{1+\beta}}([0,1])$, then 
$Y(t)$ is Markov and has a continuous version.
\end{Lemma}

\noindent 
{\em Proof:} 

     To show that a process is Markov, it is sufficient to show that its 
non-overlapping increments are independent.  For a Gaussian process 
this amounts to checking that any two 
non-overlapping increments are uncorrelated: for $s_1 < t_1 < s_2 < t_2$, 
\be
&&	E (Y(t_1) - Y(s_1)) (Y(t_2) - Y(s_2)) \\
	&& \qquad =  E Y(t_1)Y(t_2) - E Y(t_1)Y(s_2) - E Y(s_1)Y(t_2) + E Y(s_1)Y(s_2)		\\
	&& \qquad = q_f(0, t_1) - q_f(0, t_1) - q_f(0, s_1) + q_f(0, s_1) = 0.
\ee

	Next we show that $Y(t)$ has a continuous version.
\be 
&&	E(Y(t) - Y(s))^2 = q_f(0, t) - 2q_f(0, s) + q_f(0, s)		\\
	&& \qquad =C \alpha(\alpha + 1) \int_{s}^{t}\int_{s}^{t}f(u)f(v)
		\frac{du dv}{|u-v|^{1-\alpha}}+
	2 C \alpha(\alpha + 1) \int_{0}^{s}\int_{s}^{t}f(u)f(v)
		\frac{du dv}{|u-v|^{1-\alpha}}				\\
	&& \qquad = I_1 + I_2.
\ee
$I_1$ can be estimated by the Hardy-Littlewood-Sobolev 
inequality:
\begin{equation}\label{kolmogorov_estimate}
I_1 = q_f(s, t) \leq C \alpha(\alpha + 1) |s - t|^{\alpha-\beta} 
		\int_{s}^{t} \int_{s}^{t}
		\frac {f(u) f(v)} {|u-v|^{1-\beta}} du dv 
    \leq {\rm~const~}\cdot|s - t|^{\alpha-\beta} \|f\|^2_{\frac{2}{1+\alpha}}.
\end{equation}
$I_2$ can be estimated by the Hardy-Littlewood-Sobolev and H\"older 
inequalities:
\be
	I_2 &\leq& {\rm~const~}\cdot \|f \chi_{[0,s]}\|_{\frac{2}{1+\alpha}} 
		  \|f \chi_{[s,t]}\|_{\frac{2}{1+\alpha}} \\
	    &\leq& {\rm~const~}\cdot \|f \chi_{[0,s]}\|_{\frac{2}{1+\alpha}} 
		 |s - t|^{\frac{\alpha-\beta}{\alpha+\beta}}
	\bigl(\int_{s}^{t}
		  f(\tau)^{\frac{2}{1+\beta}} d\tau \bigr)^{\frac{1+\beta}
			{1+\alpha}}.
\ee
Choosing an integer $m$ such that $(\alpha-\beta)m > 1$ and 
$\frac{\alpha-\beta}{\alpha+\beta}m > 1$, we can ensure
\be
	E\left(Y(t) - Y(s)\right)^{2 m}
	\leq {\rm~const~}\cdot |s - t|^{\gamma}, 
\qquad {\rm where} \quad \gamma > 1.
\ee 
$Y$ has a continuous version by Kolmogorov's continuity criterion.
\qed

	Finally, we have a reflection principle:

\begin{Lemma}\label{reflection}

Let $Y(t)$ be a centered Gaussian Markov process with continuous
paths.  Then for all $\lambda > 0$ and $T \geq 0$,
$P(\max\limits_{0 \leq t \leq T} Y(t) \geq \lambda)
= 2 P(Y(T) \geq \lambda)$.
\end{Lemma}

\noindent
The proof is exactly analogous to the oBm case.

	Now we are ready to prove Theorem~\ref{max}.

\subsubsection{Proof of Part (1).}

\noindent 
{\em Proof:} 

Choose an integer $m$ such that $(\alpha-\beta) m > 1$.  Then
\be
	E\bigl( \int_0^t f(\tau)d\xa(\tau) -
		\int_0^s f(\tau)d\xa(\tau) \bigr)^{\! 2 m}
	&\leq& {\rm~const~}\cdot 
	\bigl( E\bigl(\int_0^t f(\tau)d\xa(\tau) -
		\int_0^s f(\tau)d\xa(\tau)\bigr)^{\! 2} \bigr)^m \\
	&=&  \big( C \alpha(\alpha + 1) \int_{s}^{t} \int_{s}^{t}
	\frac {f(u) f(v)} {|u-v|^{1-\alpha}} du dv \bigr)^m  \\
	&\leq& {\rm~const~}\cdot |s - t|^{(\alpha-\beta) m},
\ee
where the first inequality holds because
$\int_0^t f(\tau)d\xa(\tau)$ is a Gaussian random variable,
and the second by (\ref{kolmogorov_estimate}).  By Kolmogorov's 
criterion the process $\int_0^t f(\tau)d\xa(\tau)$ admits a continuous 
version.
\qed

\subsubsection{Proof of Part (2).}

\noindent 
{\em Proof:} 

Let $X(t) = \int_0^t f(\tau)d\xa(\tau)$.  $X(t)$ is a Gaussian process 
with $E X(t) = 0$ and $E X(s)X(t) = \int_{0}^{s}\int_{0}^{t}f(u)f(v)
\frac{du dv}{|u-v|^{1-\alpha}}$. 
Define $Y(t)$ to be a Gaussian process with $E Y(t) = 0$ and 
$E Y(s)Y(t) = q_f(0,s)$ for $s \leq t$.  Clearly $E X(t)^2 = E Y(t)^2$. 
It can be shown that the process $Y(t)$ is well-defined for all $t$. 

Suppose $f \geq 0$.  Then $E X(s)X(t) \geq E Y(s)Y(t)$ for 
$s, t \in [0,1]$.
Therefore the processes $X(t)$ and $Y(t)$ satisfy the assumptions of 
Slepian's lemma (Lemma~\ref{Slepian}). 
By Lemma~\ref{Ymarkov-continuous}, $Y(t)$ is a Markov process with 
continuous paths, and by Lemma~\ref{reflection},
$Y$ obeys the reflection principle:
\be
  P(\max\limits_{0 \leq t \leq 1} Y(t) \geq \lambda) =
  2 P(Y(1) \geq \lambda) 	=
\int_{\lambda/\sqrt{q_f(0,1)}}^{\infty} \sqrt{\frac{2}{\pi}}\, e^{-x^2/2} dx.
\ee
Hence for $f \geq 0$,
\be
   P(\max\limits_{0 \leq t \leq 1} X(t) \geq \lambda)
   \leq 
\int_{\lambda/\sqrt{q_f(0,1)}}^{\infty} \sqrt{\frac{2}{\pi}}\, e^{-x^2/2} dx.
\ee
Let $f \in L^{\frac{2}{1+\beta}}([0, \infty))$.  
Write $X_{\pm} = \int_{0}^{t} f_{\pm}(\tau) d\xa(\tau)$.
Define processes $Y_{\pm}(t)$ by replacing $f$ by $f_{\pm}$ in the 
definition of $Y(t)$.  Since Slepian's lemma applies
also to $-X_{-}(t)$ and $Y_{-}(t)$, we have
\be
	P(\max\limits_{0 \leq t \leq 1} \pm X_{\pm}(t) \geq \lambda)
	\leq 
 \int_{\lambda/\sqrt{q_{f_{\pm}}(0,1)}}^{\infty} \sqrt{\frac{2}{\pi}} \,
						e^{-x^2/2} dx.
\ee
Therefore
\be
 P(\max\limits_{0 \leq t \leq 1} X(t) \geq \lambda) &=&
	P(\max\limits_{0 \leq t \leq 1} X_{+}(t) 
	+ \max\limits_{0 \leq t \leq 1} (- X_{-}(t)) \geq \lambda) \\
      &\leq&
	P(\max\limits_{0 \leq t \leq 1} X_{+}(t) \geq \lambda r) +
	P(\max\limits_{0 \leq t \leq 1} (-X_{-}(t)) \geq \lambda(1 - r)) \\    
&\leq& 	\int_{\lambda r/\sqrt{q_{f_+}(0,1)}}^{\infty} +
		\int_{\lambda (1 - r)/\sqrt{q_{f_-}(0,1)}}^{\infty}
				\sqrt{\frac{2}{\pi}}\, e^{-x^2/2} dx.
\ee
\qed

\subsubsection{Proof of Part (3).}

\noindent 
{\em Proof:} 

Inequality~(\ref{kolmogorov_estimate}) shows that
$\sqrt{q_f(s, t)} \leq {\rm~const~}\cdot |s - t|^{\frac{\alpha-\beta}{2}}$,
so $$\int_1^{\infty} \sup\limits_{|s-t| 
	< m^{-x^2}} \sqrt{q_f(s, t)}\, dx < \infty.$$ 
This is the condition of applicability of Fernique's inequality
\cite{F}, of which the claim is a direct consequence.
\qed


\begin{thebibliography}{XX}
\bibitem[Dy]{Dy} A.M. Dyachkov, Conditions for the existence of Stieltjes integral of functions of bounded generalized variation, {\it Anal. Math.} {\bf 14}
(1988) no.4 295--313.

\bibitem[DH]{DH} W. Dai and C. C. Heyde, It\^o's formula with respect 
to fractional Brownian motion and its application, 
{\it J. Appl. Math. Stochastic Anal.} {\bf 9} (1996) 439--448. 

\bibitem[DU]{DU} L. Decreusefond and  \"Ust\"unel, Ali S\"uleyman,
Application du calcul des variations stochastiques au mouvement brownien 
fractionnaire, {\it C. R. Acad. Sci. Paris}, s\'er.\ I {\bf 321} (1995)
1605--1608.

\bibitem[F]{F} X. Fernique, R\'egularit\'e des trajectoires des fonctions 
al\'eatoires gaussiennes, {\it Ecole d'Et\'e de Probabilit\'es de Saint-Flour, 
avril 1974}, Lecture Notes in Math.\ vol.\ 480, Springer, Berlin, 
1975, 47--52. 

\bibitem[G]{G} A. Goldman, Estimations analytiques concernant le mouvement
brownien fractionnaire \`a plusieurs param\`etres {\it
C. R. Acad. Sci. Paris}, s\'er.\ I {\bf 298} (1984) 91--93. 

\bibitem[K1]{K1} J.-P. Kahane, {\it Some random series of functions},  
Cambridge Studies in Adv. Math.\ 5, 1985.

\bibitem[K2]{K2} J.-P. Kahane,  Sur les mouvements browniens 
fractionnaires: images, graphes, niveaux, {\it C. R.
Acad. Sci. Paris}, s\'er.\ I {\bf  300} (1985) 501--503. 

\bibitem[Kl]{Kl} A. N. Kolmogorov, Wiener screw-line and other 
interesting curves in Hilbert space, {\it Doklady Akad. Nauk\/}
{\bf 26} (1940) 115--118.

\bibitem[Kn]{Kn} V. Kondurar, Sur l'int\'egrale de Stieltjes, 
{\it Mat. Sbornik\/} {\bf 2} (1937) 361--366. 

\bibitem[Lb]{Lb} A. Le Breton, Filtering and parameter estimation 
in a simple linear system driven by a fractional Brownian motion, 
{\it Statist. Probab. Lett.} {\bf 38} (1998) 263--274.  

\bibitem[LL]{LL} E.H. Lieb and M. Loss, {\it Analysis},  Graduate Studies in Mathematics 14, {\it AMS}, 1997.

\bibitem[Ln]{Ln} S. J. Lin, Stochastic analysis of fractional Brownian 
motions, {\it Stochastics Rep.} {\bf 55} (1995) 121--140.

\bibitem[Lv]{Lv} P. L\'evy, {\it Processus stochastiques et movement 
brownien}, Paris, 1948 (2nd ed.\ 1965).
 
\bibitem[Mk]{Mk} H. P. McKean Jr, {\it Stochastic integrals}, Probability and Mathematical Statistics 5, {\it Academic Press}, 1969.

\bibitem[Mn]{Mn} B. Mandelbrot and J. W. Van Ness, Fractional Brownian 
motions, fractional noises and applications, {\it SIAM Rev.} {\bf 10} 
(1968) 422--437. 

\bibitem[Ml]{Ml} G. M. Molchan, Maximum of fractional Brownian 
motion: probabilities of small values (preprint).

\bibitem[Ru]{Ru} A.A. Ruzmaikina, Stochastic Calculus with fractional Brownian motion, Ph.D. Thesis, Princeton University (1999).

\bibitem[Sn]{Sn} Ya. G. Sina\u\i, On the distribution of the maximum 
of fractional Brownian motion, {\it Uspekhi Mat. Nauk\/} {\bf 52} (1997)
119--138.

\bibitem[Sl]{Sl} D. Slepian, The one-sided barrier problem for Gaussian 
noise, {\it Bell System Tech. J.} {\bf 41} (1962) 463--501. 

\bibitem[T1]{T1} M. Talagrand, Hausdorff measure of trajectories of 
multiparameter fractional Brownian motion, {\it Ann. Probab.} 
{\bf 23} (1995) 767--775.

\bibitem[T2]{T2} M. Talagrand, Lower classes for fractional Brownian 
motion, {\it J. Theoret. Probab.} {\bf 9} (1996) 191--213.  

\bibitem[W]{W} M. Weber, Dimension de Hausdorff et points multiples 
du mouvement brownien fractionnaire dans $\R^n$, {\it
C. R. Acad. Sci. Paris}, s\'er.\ I {\bf 297} (1983) 357--360.

\bibitem[X]{X} Y. Xiao, Hausdorff-type measures of the sample paths 
of fractional Brownian motion, {\it Stochastic Process. Appl.} 
{\bf 74} (1998) 251--272.

\bibitem[Yr]{Yr} M. Yor, Remarques sur certaines constructions des 
mouvements browniens fractionnaires, {\it S\'eminaire de
Probabilit\'es, XXII}, Lecture Notes in Math.\ vol.\ 1321, Springer,
Berlin-New York, 1988, 217--224.

\bibitem[Yg]{Yg} L. C. Young, An inequality of the H\"older type connected 
with Stieltjes integration, {\it Acta Math.} {\bf 67} (1936) 251--282.

\end{thebibliography}
 \end{document}